\documentclass[12pt]{article}
\usepackage{epic,latexsym}
\usepackage{color}
\usepackage{tikz}
\usepackage{comment}
\usepackage{amssymb}
\usepackage{mathtools}
\usepackage{url}

\newtheorem{theorem}{Theorem}[section]

\newtheorem{proposition}[theorem]{Proposition}

\newtheorem{corollary}[theorem]{Corollary}

\newtheorem{lemma}[theorem]{Lemma}

\newcommand{\gp}{{\rm gp}}
\newcommand{\gps}{{\rm gp}{_{\rm S}}}
\newcommand{\gpsl}{\underline{{\rm gp}}{_{\rm S}}}
\newcommand{\SP}{Sierpi\'nski product}
\newcommand{\proof}{\noindent{\bf Proof.\ }}
\newcommand{\qed}{\hfill $\square$ \bigskip}
\newcommand{\cp}{\,\square\,}

\date{}

\title{General position sets, colinear sets, and Sierpi\'{n}ski product graphs}

\author{
Jing Tian $^{a,b}$ \and Sandi Klav\v{z}ar $^{b,c,d}$\\\\
$^{a}$ \small School of Science, Zhejiang University of Science and Technology, \\
\small Hangzhou, Zhejiang 310023, PR China\\
\small {\tt jingtian526@126.com}\\
$^{b}$ \small Institute of Mathematics, Physics and Mechanics, Ljubljana, Slovenia \\
$^{c}$\small Faculty of Mathematics and Physics, University of Ljubljana, Slovenia\\
$^{d}$ \small Faculty of Natural Sciences and Mathematics, University of Maribor, Slovenia\\
\small {\tt sandi.klavzar@fmf.uni-lj.si}\\
}

\begin{document}

\maketitle

\begin{abstract}
Let $G \otimes _f H$ denote the Sierpi\'nski product of graphs $G$ and $H$ with respect to the function $f$. The Sierpi\'nski general position number ${\rm gp}{_{\rm S}}(G,H)$ is introduced as the cardinality of a largest general position set in $G \otimes _f H$ over all possible functions $f$. Similarly, the lower Sierpi\'nski general position number $\underline{{\rm gp}}{_{\rm S}}(G,H)$ is the corresponding smallest cardinality. The concept of vertex-colinear sets is introduced. Bounds for the general position number in terms of extremal vertex-colinear sets, and bounds for the (lower) Sierpi\'nski general position number are proved. The extremal graphs are investigated. Formulas for the (lower) Sierpi\'nski general position number of the \SP{s} with $K_2$ as the first factor are deduced. It is proved that if $m,n\geq 2$, then ${\rm gp}{_{\rm S}}(K_m,K_n) = m(n-1)$ and that if $n\ge 2m-2$, then $\underline{{\rm gp}}{_{\rm S}}(K_m,K_n) = m(n-m+1)$.
\end{abstract}


\medskip\noindent
{\bf AMS Subj.\ Class.\ (2020)}: 05C12, 05C76

%%%%%%%%%%%%%%%%%%%%%%
\section{Introduction}
%%%%%%%%%%%%%%%%%%%%%%

General position sets were introduced to graph theory in~\cite{manuel-2018}, but Chandran S.V.\ and Parthasarathy in~\cite{ullas} earlier studied the concept of geodetic irredundant sets which turned out to be an equivalent concept to general position sets. Moreover, much earlier K\"orner~\cite{Korner-1995} studied the general position sets of hypercubes (although in a completely different context). Briefly, a general position set is a set of vertices such that no three vertices from the set lie on a common shortest path.

The problem of finding a cardinality of a largest general position set in a graph, that is, the general position number of the graph, is NP-hard~\cite{manuel-2018}. In~\cite{anand-2019}, general position sets of a graph were characterized. Afterwards, the problem of determining the general position number of (classes of) graphs received extensive attention, see~\cite{Irsic-2023, Patkos-2019, tian-2023, yao-2022}. General position sets were generalized to $d$-position sets~\cite{klavzar-2021} and to Steiner general position sets~\cite{klavzar-2021a}. In addition, in~\cite{stefano-2024+} the lower general position was investigated, that is, the cardinality of a smallest maximal general position set in a graph. The edge version of general position sets were investigated in~\cite{Klavzar-Tan-2023, manuel-2022, tian-2024}, while for the  monophonic version see~\cite{thomas-2024+}.

A lot of attention has been given to the general position sets in Cartesian  product graphs, see~\cite{Klavzar-Rus-2021, Korze-Vesel-2023, Tian-2021}. In this paper we are interested in general position sets in the recently introduced graph operation called the Sierpi\'{n}ski product. The operation was introduced by Kovi\v{c}, Pisanski, Zemlji\v{c}, and \v{Z}itnik in~\cite{kpzz-2022} with the idea to generalize the Sierpi\'{n}ski graphs~\cite{klavzar-1997}, cf.~Fig.~\ref{fig:K4-otimes-K4}. The investigation of the latter graphs up to 2017 is summarized in the survey~\cite{hinz-2017}, some of the recent research of Sierpi\'{n}ski graphs can be found, for instance, in~\cite{menon-2023, palma-2024, varghese-2024}, see also the references therein. The Sierpi\'{n}ski product was further investigated in~\cite{henning-2024+, henning-2024}, where the (upper) Sierpi\'{n}ski domination number and the (upper) Sierpi\'{n}ski metric dimension were respectively introduced. Following this line of investigation, in this paper we introduce the Sierpi\'nski general position number $\gps(G,H)$, and the lower Sierpi\'nski general position number $\gpsl(G,H)$, of the Sierpi\'nski product of graphs $G$ and $H$.

In the rest of the paper we proceed as follows. In the next section we give definitions and recall known results that we will use in the paper. In Section~\ref{sec:colinear} we introduce a new concept into the theory of graph general position sets, the vertex-colinear sets. This concept is useful while studying the general position sets in \SP\ graphs, but we believe to be also of independent interest. In particular, we prove bounds for the general position number of a graph in terms of extremal vertex-colinear sets and show their usefulness for graphs with bridges. In Section~\ref{sec:arbitrary}, we present bounds for the (lower) Sierpi\'nski general position number for general graphs and give formulas for these numbers for the \SP\ graphs with $K_2$ as the first factor. In Section~\ref{sec:complete graphs}, we consider the \SP\ of complete graphs $K_m \otimes _f K_n$. We determine their Sierpi\'nski general position number, more precisely, if $m,n\geq 2$, then $\gps(K_m,K_n) = m(n-1)$. For the lower Sierpi\'nski general position number we prove that $\gpsl(K_m,K_n) = m(n-m+1)$ provided that $n\ge 2m-2$. The  assumption $n\ge 2m-2$ is required as demonstrated by our last result which asserts that $\gpsl(K_6,K_9)=25$.

%%%%%%%%%%%%%%%%%%%%%%
\section{Preliminaries}
\label{sec:preliminaries}
%%%%%%%%%%%%%%%%%%%%%%

In this section we define concepts, notation, and results needed, and introduce the (lower) Sierpi\'{n}ski general position number. We begin with basic definitions, and follow by introducing the general position sets and the \SP\ graphs.

For a positive integer $k$ we set $[k]=\{1,\ldots,k\}$. Unless stated otherwise, graphs $G=(V(G),E(G))$ in the paper are connected and simple. The {\em order} $n(G)$ of $G$ is equal to $|V(G)|$. The {\em degree} of a vertex $u$, $\deg_G(u)$,
is the number of vertices adjacent to $u$ in $G$. Vertices of degree one are called {\em leaves}. The number of leaves of $G$ will be denoted by $\ell(G)$. For a graph $G$, let $v_k(G)$ denote the number of vertices of $G$ of degree $k$, cf.~\cite{cai-1992}. An edge $e$ of $G$ is a {\em bridge} if $G-e$ is disconnected.  If $S\subseteq V(G)$, then the subgraph of $G$ induced by $S$ is denoted by $G[S]$. A vertex of $G$ is {\em simplicial} if its neighbourhood induces a complete subgraph. The {\em distance} $d_G(u,v)$ between vertices $u$ and $v$ of $G$ is the number of edges on a shortest $u,v$-path. The {\em interval} between vertices $u$ and $v$ is
$$I_G[u,v] = \{w:\ d_G(u,v) = d_G(u,w) + d_G(w,v)\}\,.$$
A subgraph $H$ of $G$ is {\em isometric} if for each pair of vertices $u,v\in V(H)$ we have $d_H(u,v) = d_G(u,v)$ and is {\em convex} if whenever $u,v\in V(H)$ and $P$ is a shortest $u,v$-path in $G$, then $P$ lies completely in $H$.

Let $X\subseteq V(G)$. Then $X$ is a {\em general position set} of $G$ if for any shortest path $P$ in $G$ we have $|V(P)\cap X| \le 2$. The cardinality of a largest general position set of $G$ is the {\em general position number}, $\gp(G)$, of $G$. A general position set $X$ of cardinality $\gp(G)$ is referred to as a {\em $\gp$-set} of $G$. Clearly, if $G$ has at least two vertices, then $\gp(G) \ge 2$.

Let $G$ and $H$ be graphs and let $f \colon V(G)\rightarrow V(H)$ be a function. The {\em Sierpi\'{n}ski product} of $G$ and $H$ (with respect to $f$) is the graph  $G \otimes _f H$ with vertices
\begin{itemize}
 \item $V(G \otimes _f H) = V(G)\times V(H)$,
\end{itemize}
and with edges
\begin{itemize}
    \item $(g,h)(g,h')$, where $g\in V(G)$ and $hh' \in E(H)$, and
    \item $(g,f(g'))(g',f(g))$, where $gg' \in E(G)$.
\end{itemize}
The second type edges that we have just defined  will be called the \emph{connecting edges} of $G \otimes _f H$. Note that for each vertex $g \in V(G)$, the subgraph of $G \otimes _f H$ induced by the set of vertices $\{(g,h):\ h\in V(H)\}$, is isomorphic to $H$; it will be denoted by $gH$.

For an example consider the Sierpi\'{n}ski product of $G\cong K_4$ and $H\cong K_4$, where $V(G) = [4]$, $V(H) = \{x_1, x_2, x_3, x_4\}$ and $f:V(G)\rightarrow V(H)$ is defined with $f(i) = x_i$, $i\in [4]$. On the left side of Fig.~\ref{fig:K4-otimes-K4} the graph $G \otimes _f H$ is drawn according to the graph definition, while on the right side the standard image of the Sierpi\'{n}ski graph $S_4^2$ (cf.~\cite{hinz-2017}) is presented. 

\begin{figure}[ht!]
\begin{center}
\begin{tikzpicture}[scale=0.7,style=thick,x=1cm,y=1cm]
\def\vr{3pt}

\begin{scope}[xshift=0cm, yshift=0cm] % K4 1
\coordinate(x1) at (0,0);
\coordinate(x2) at (0,2);
\coordinate(x3) at (0,4);
\coordinate(x4) at (0,6);
% \edges		
\draw (x1) -- (x2) -- (x3) -- (x4);
\draw (x1) .. controls (0.4,1) and (0.4,2) .. (x3);
\draw (x2) .. controls (-0.4,3) and (-0.4,4) .. (x4);
\draw (x1) .. controls (0.8,1) and (0.8,5) .. (x4);
\draw (0,2) -- (2,0);
\draw (2,4) -- (4,2);
\draw (4,6) -- (6,4);
\draw (0,4) .. controls (0.5,1.5) and (3,0.5) .. (4,0);
\draw (2,6) .. controls (3,5.5) and (5.5,3.5) .. (6,2);
\draw (0,6) .. controls (1,2) and (5,1) .. (6,0);
%  vertices
\draw(x1)[fill=white] circle(\vr);
\draw(x2)[fill=white] circle(\vr);
\draw(x3)[fill=white] circle(\vr);
\draw(x4)[fill=white] circle(\vr);
% text
\node at (0,-1) {$1$};
\node at (2,-1) {$2$};
\node at (4,-1) {$3$};
\node at (6,-1) {$4$};
\node at (-1,0) {$x_1$};
\node at (-1,2) {$x_2$};
\node at (-1,4) {$x_3$};
\node at (-1,6) {$x_4$};
\node at (8,3) {$\cong$};
\end{scope}
\begin{scope}[xshift=2cm, yshift=0cm] % K4 2
\coordinate(x1) at (0,0);
\coordinate(x2) at (0,2);
\coordinate(x3) at (0,4);
\coordinate(x4) at (0,6);
% \edges		
\draw (x1) -- (x2) -- (x3) -- (x4);
\draw (x1) .. controls (0.4,1) and (0.4,2) .. (x3);
\draw (x2) .. controls (-0.4,3) and (-0.4,4) .. (x4);
\draw (x1) .. controls (0.8,1) and (0.8,5) .. (x4);
%  vertices
\draw(x1)[fill=white] circle(\vr);
\draw(x2)[fill=white] circle(\vr);
\draw(x3)[fill=white] circle(\vr);
\draw(x4)[fill=white] circle(\vr);
% text
% \node at (0.2,0.75) {$u$};
\end{scope}
\begin{scope}[xshift=4cm, yshift=0cm] % K4 3
\coordinate(x1) at (0,0);
\coordinate(x2) at (0,2);
\coordinate(x3) at (0,4);
\coordinate(x4) at (0,6);
% \edges		
\draw (x1) -- (x2) -- (x3) -- (x4);
\draw (x1) .. controls (0.4,1) and (0.4,2) .. (x3);
\draw (x2) .. controls (-0.4,3) and (-0.4,4) .. (x4);
\draw (x1) .. controls (0.8,1) and (0.8,5) .. (x4);
%  vertices
\draw(x1)[fill=white] circle(\vr);
\draw(x2)[fill=white] circle(\vr);
\draw(x3)[fill=white] circle(\vr);
\draw(x4)[fill=white] circle(\vr);
% text
% \node at (0.2,0.75) {$u$};
\end{scope}
\begin{scope}[xshift=6cm, yshift=0cm] % K4 4
\coordinate(x1) at (0,0);
\coordinate(x2) at (0,2);
\coordinate(x3) at (0,4);
\coordinate(x4) at (0,6);
% \edges		
\draw (x1) -- (x2) -- (x3) -- (x4);
\draw (x1) .. controls (0.4,1) and (0.4,2) .. (x3);
\draw (x2) .. controls (-0.4,3) and (-0.4,4) .. (x4);
\draw (x1) .. controls (0.8,1) and (0.8,5) .. (x4);
%  vertices
\draw(x1)[fill=white] circle(\vr);
\draw(x2)[fill=white] circle(\vr);
\draw(x3)[fill=white] circle(\vr);
\draw(x4)[fill=white] circle(\vr);
% text
% \node at (0.2,0.75) {$u$};
\end{scope}
%
% SECOND DRAWING
%
\begin{scope}[xshift=10cm, yshift=0cm] % K4 1
\coordinate(x1) at (0,0);
\coordinate(x2) at (0,2);
\coordinate(x3) at (2,2);
\coordinate(x4) at (2,0);
% \edges		
\draw (x1) -- (x2) -- (x3) -- (x4) -- (x1) -- (x3);
\draw (x2) -- (x4);
\draw (2,0) -- (4,0);
\draw (2,6) -- (4,6);
\draw (0,2) -- (0,4);
\draw (6,2) -- (6,4);
\draw (2,2) -- (4,4);
\draw (2,4) -- (4,2);
%  vertices
\draw(x1)[fill=white] circle(\vr);
\draw(x2)[fill=white] circle(\vr);
\draw(x3)[fill=white] circle(\vr);
\draw(x4)[fill=white] circle(\vr);
\end{scope}
\begin{scope}[xshift=14cm, yshift=0cm] % K4 2
\coordinate(x1) at (0,0);
\coordinate(x2) at (0,2);
\coordinate(x3) at (2,2);
\coordinate(x4) at (2,0);
% \edges		
\draw (x1) -- (x2) -- (x3) -- (x4) -- (x1) -- (x3);
\draw (x2) -- (x4);
%  vertices
\draw(x1)[fill=white] circle(\vr);
\draw(x2)[fill=white] circle(\vr);
\draw(x3)[fill=white] circle(\vr);
\draw(x4)[fill=white] circle(\vr);
\end{scope}
\begin{scope}[xshift=10cm, yshift=4cm] % K4 3
\coordinate(x1) at (0,0);
\coordinate(x2) at (0,2);
\coordinate(x3) at (2,2);
\coordinate(x4) at (2,0);
% \edges		
\draw (x1) -- (x2) -- (x3) -- (x4) -- (x1) -- (x3);
\draw (x2) -- (x4);
%  vertices
\draw(x1)[fill=white] circle(\vr);
\draw(x2)[fill=white] circle(\vr);
\draw(x3)[fill=white] circle(\vr);
\draw(x4)[fill=white] circle(\vr);
\end{scope}
\begin{scope}[xshift=14cm, yshift=4cm] % K4 4
\coordinate(x1) at (0,0);
\coordinate(x2) at (0,2);
\coordinate(x3) at (2,2);
\coordinate(x4) at (2,0);
% \edges		
\draw (x1) -- (x2) -- (x3) -- (x4) -- (x1) -- (x3);
\draw (x2) -- (x4);
%  vertices
\draw(x1)[fill=white] circle(\vr);
\draw(x2)[fill=white] circle(\vr);
\draw(x3)[fill=white] circle(\vr);
\draw(x4)[fill=white] circle(\vr);
\end{scope}

\end{tikzpicture}
\caption{The Sierpi\'{n}ski product $K_4 \otimes _f K_4$, where $f(i) = x_i$, $i\in [4]$.}
	\label{fig:K4-otimes-K4}
\end{center}
\end{figure}

Let $G$ and $H$ be connected graphs and denote by $H^G$ the set of functions from $V(G)$ to $V(H)$. We introduce the \textit{Sierpi\'{n}ski general position number}, $\gps(G,H)$, as the cardinality of a largest general position set in $G \otimes _f H$ over all possible functions $f\in H^G$ and the {\em lower Sierpi\'nski general position number}, $\gpsl(G,H)$, as the corresponding smallest cardinality. That is,
\[ \gps(G, H) = \max_{f\in H^G}\{\gp(G\otimes _f H)\} \qquad
{\rm and} \qquad
\gpsl(G, H) = \min_{f\in H^G}\{\gp(G\otimes _f H)\}\,.
\]

At the end of the section, we will recall several known results needed later on. For the first one, known as the {\em Isometric Cover Lemma}, we need the following definition. A set of subgraphs $\{H_1,\ldots, H_k\}$ is an {\em isometric cover} of $G$ if each $H_i$, $i\in [k]$, is isometric in $G$ and $\bigcup_{i=1}^k V(H_i) = V(G)$.

\begin{theorem} {\rm \cite[Theorem 3.1]{manuel-2018}}
\label{thm:iso-cover}
If $\{H_1,\ldots, H_k\}$ is an isometric cover of $G$, then
$$\gp(G) \le \sum_{i=1}^k \gp(H_i)\,.$$
\end{theorem}

\begin{theorem} {\rm \cite[Theorem 3.6]{manuel-2018}}
\label{thm:block-graphs}
If $S$ is the set of simplicial vertices of a block graph $G$, then $\gp(G) = |S|$.
\end{theorem}

If $G$ is a connected graph, $S\subseteq V(G)$, and ${\cal P} = \{S_1, \ldots, S_t\}$ a partition of $S$, then ${\cal P}$ is \emph{distance-constant} if for any $i,j\in [t]$, $i\ne j$, the distance $d_G(x,y)$, where $x\in S_i$ and $y\in S_j$, is independent of the selection of $x$ and $y$. A distance-constant partition ${\cal P}$ is {\em in-transitive} if $d_G(S_i, S_k) \ne d_G(S_i, S_j) + d_G(S_j,S_k)$ holds for $i,j,k\in [t]$. The following characterization of general position sets will be used either implicitly or explicitly in the rest of the paper.

\begin{theorem} {\rm \cite[Theorem 3.1]{anand-2019}}
\label{thm:gpsets}
Let $G$ be a connected graph. Then $S\subseteq V(G)$ is a general position set if and only if the components of $G[S]$ are complete subgraphs, the vertices of which form an in-transitive, distance-constant partition of $S$.
\end{theorem}

The last result needed concerns the metric structure of Sierpi\'{n}ski product graphs.

\begin{theorem} {\rm \cite[Theorem 4.1]{henning-2024}}
\label{thm:convex}
If $G$ and $H$ are connected graphs, $f \in H^G$, and $g\in V(G)$, then $gH$ is a convex subgraph of $G \otimes _f H$.
\end{theorem}

%%%%%%%%%%%%%%%%%%%%%%
\section{Colinear sets}
\label{sec:colinear}
%%%%%%%%%%%%%%%%%%%%%%

In this section we introduce a new concept into the theory of graph general position sets, the vertex-colinear sets. We believe that the concept is interesting in its own right. Moreover, from our perspective it turns out to be extremely useful for studying the general position sets in \SP\ graphs.

Let $G$ be a connected graph and $u\in V(G)$. A set $S\subseteq V(G)$ is a {\em $u$-colinear set} if $S$ is a general position set such that $u\notin S$ and $y\notin I_G[x,u]$ for any $x,y\in S$. Before moving on, we point out that a similar concept named {\em $u$-position set} was independently introduced in~\cite{thankachy-2024}. More precisely, a set $S$ is a $u$-colinear set if and only if $S$ is both a general position set and a $u$-position set.

We are interested in $u$-colinear sets of largest cardinality and hence set
$$\xi_G(u) = \max \{ |S|:\ S\ {\rm is\ a}\ u\mbox{-}{\rm colinear\ set}\}\,.$$
Further, let
\begin{align*}
\xi^{-}(G) & = \min \{\xi_G(u):\ u\in V(G)\}\,,\\
\xi(G) & = \max \{\xi_G(u):\ u\in V(G)\}\,.
\end{align*}

To illustrate these concepts, we state the following easy result which follows from \cite[Corollary~19]{thankachy-2024}. To be self-contained, we include its proof here.

\begin{proposition}
\label{prop:xi-in-tree}
If $T$ is a tree with $n(T) \ge 3$ and $u\in V(T)$, then $\xi_T(u) = \ell(T) - 1$ if $u$ is a leaf, and $\xi_T(u) = \ell(T)$, otherwise. In particular, $\xi^{-}(T) = \ell(T) - 1$ and $\xi(T) = \ell(T)$.
\end{proposition}

\proof
Assume first that $u$ is a leaf and let $S$ be a $u$-colinear set. Let $u'$ be the support vertex adjacent to $u$. In the following arguments we use the the consequence of Theorem~\ref{thm:block-graphs} asserting that the general position number of a tree is the number of its leaves. If $\deg_T(u') > 2$, then $\gp(T-u) = \ell(T) - 1$, and then the leaves of $T-u$ form a largest $u$-colinear set. Assume next that $\deg_T(u') = 2$. Then $T-u$ contains $\ell(T)$ leaves (one of them being $u'$). However, if $u'$ lies in a $u$-colinear set, then no other vertex lies in such a set. Hence we can again conclude that $\xi_T(u) = \ell(T) - 1$. Since $n(T)\ge 3$, there exists a vertex $u$ which is not a leaf. Then the leaves of $T$ form a $u$-colinear set and hence $\xi_T(u) = \ell(T)$ and so also $\xi(T) = \ell(T)$.
\qed

The above invariants fulfil the following chain of inequalities involving the general position number. This in part demonstrates that our definitions are meaningful, further reasons for their introduction will be given later on.

\begin{theorem}
\label{thm:2xi>gp}
If $G$ is a connected graph of order at least $2$ and $u\in V(G)$, then
$$\xi^{-}(G) \le \xi_G(u) \le \xi(G)\le \gp(G) \leq 2\xi^{-}(G)\,.$$
\end{theorem}

\proof
The left three inequalities follow directly from the definition. To prove the last inequality, let $u$ be a vertex of $G$ such that $\xi_G(u) = \xi^{-}(G)$ and let $X$ be an arbitrary gp-set of $G$. We distinguish two cases.

Assume first that $u\in X$. Since $X$ is a general position set of $G$, we infer that $X \setminus \{u\}$ is a $u$-colinear set. Then it follows that $\xi_G(u)\geq |X|-1$, hence $|X|\leq \xi^{-}(G)+1$. Since the order of $G$ is at least two,  $\xi^{-}(G)\geq 1$ and thus we have $|X|\leq \xi^{-}(G)+\xi^{-}(G)\leq 2\xi^{-}(G)$.

Assume second that $u\notin X$. In this case we partition $X$ into two sets $X_1$ and $X_2$ as follows. For every $x\in X$, we put each vertex from $(I_G[x,u]\setminus \{x\})\cap X$ into $X_2$. This defines $X_2$, and then $X_1 = X\setminus X_2$. We claim that each of $X_1$ and $X_2$ is a $u$-colinear set. By the definition, this clearly holds for $X_1$. Consider now two vertices $x,x'\in X_2$ and suppose by way of contradiction that there exists a shortest $x,u$-path $P_{x,u}$ that contains $x'$. Then $d_G(x,u) = d_G(x,x') + d_G(x',u)$. Since $x\in X_2$, there exists a vertex $x''\in X$ and a shortest $x'',u$-path $P_{x'',u}$. We are going to show that $x''$, $x$, and $x'$ lie on a common shortest path. If not, then $d_G(x'',x') < d_G(x'',x) + d_G(x,x')$, but then
\begin{align*}
d_G(x'',u) & = d_G(x'',x) + d_G(x,u) \\
           & = d_G(x'',x) + d_G(x,x') + d_G(x',u) \\
           & > d_G(x'',x') + d_G(x',u) \\
           & = d_G(x'',u)\,,
\end{align*}
which is not possible. We have thus seen that also $X_2$ is a $u$-colinear set. Therefore,
$$\xi_G(u)\geq \max\{|X_1|, |X_2|\} \geq \frac{1}{2} |X| = \frac{1}{2}\gp(G)\,,$$
hence $\gp(G) \le 2\xi_G(u) = 2\xi^{-}(G)$.
\qed

We remark that the first part of the proof of Theorem~\ref{thm:2xi>gp} can also be deduced from the proof of~\cite[Lemma~2]{thankachy-2024}.
	
We now give some examples which demonstrate that in the inequality chain of Theorem~\ref{thm:2xi>gp} the inequalities can be sharp or not. If $n\ge 2$, then $\xi^{-}(K_{n,n}) = \xi(K_{n,n}) = \gp(K_{n,n})$, hence for these graphs the left three equalities hold. On the other hand, $\xi(G)$ can be smaller than $\gp(G)$. For instance, $\xi(K_n) = n-1 = \gp(K_n)-1$ for $n\ge 2$. In addition, if $P$ is the Petersen graph, then $\xi^{-}(P) = \xi(P) = 5$ while $\gp(P) = 6$.

Consider next the cycle chain graphs $C_{2k}^\ell$, $k, \ell\ge 2$, where  $C_{2k}^\ell$ consists of $\ell$ cycles $C_{2k}$ sharing a vertex, such that a middle cycle shares its diametral vertices with its two neighboring cycles. See Fig.~\ref{fig:C_{6}^5} where $C_{6}^5$ is shown. Consider the vertex $u$ of $C_{2k}^\ell$ as shown in Fig.~\ref{fig:C_{6}^5}. Then $\xi_{C_{2k}^\ell}(u) = 2$ and hence $\xi^{-}(C_{2k}^\ell) = 2$. On the other hand, for the vertex $v$ we have $\xi_{C_{2k}^\ell}(v) = 4$ which in turn implies that $\xi(C_{2k}^\ell) = 4$. Thus $\xi^{-}(G)$ can be strictly smaller than $\xi(G)$. Moreover, the vertex $v$ of $C_{2k}^\ell$ also demonstrates that in Theorem~\ref{thm:2xi>gp}, $\xi_G(u)$ can be equal to $2\xi^{-}(G)$.

\begin{figure}[ht!]
\begin{center}
\begin{tikzpicture}[scale=1.0,style=thick,x=1cm,y=1cm]
\def\vr{3pt}

\begin{scope}[xshift=0cm, yshift=0cm] % C4 1
\coordinate(x1) at (0.5,0.75);
\coordinate(x2) at (1,0);
\coordinate(x3) at (2,0);
\coordinate(x4) at (2.5,0.75);
\coordinate(x5) at (2,1.5);
\coordinate(x6) at (1,1.5);
% \edges		
\draw (x1) -- (x2) -- (x3) -- (x4) --(x5)-- (x6) --(x1);
%  vertices
\draw(x1)[fill=white] circle(\vr);
\draw(x2)[fill=white] circle(\vr);
\draw(x3)[fill=white] circle(\vr);
\draw(x4)[fill=white] circle(\vr);
\draw(x5)[fill=white] circle(\vr);
\draw(x6)[fill=white] circle(\vr);
% text
\node at (0.2,0.75) {$u$};
\node at (2.5,0.4) {$v$};
\end{scope}

\begin{scope}[xshift=2cm, yshift=0cm] % C4 2
\coordinate(x1) at (0.5,0.75);
\coordinate(x2) at (1,0);
\coordinate(x3) at (2,0);
\coordinate(x4) at (2.5,0.75);
\coordinate(x5) at (2,1.5);
\coordinate(x6) at (1,1.5);
% \edges		
\draw (x1) -- (x2) -- (x3) -- (x4) --(x5)-- (x6) --(x1);
%  vertices
\draw(x1)[fill=white] circle(\vr);
\draw(x2)[fill=white] circle(\vr);
\draw(x3)[fill=white] circle(\vr);
\draw(x4)[fill=white] circle(\vr);
\draw(x5)[fill=white] circle(\vr);
\draw(x6)[fill=white] circle(\vr);
\end{scope}

\begin{scope}[xshift=4cm, yshift=0cm] % C4 3
\coordinate(x1) at (0.5,0.75);
\coordinate(x2) at (1,0);
\coordinate(x3) at (2,0);
\coordinate(x4) at (2.5,0.75);
\coordinate(x5) at (2,1.5);
\coordinate(x6) at (1,1.5);
% \edges		
\draw (x1) -- (x2) -- (x3) -- (x4) --(x5)-- (x6) --(x1);
%  vertices
\draw(x1)[fill=white] circle(\vr);
\draw(x2)[fill=white] circle(\vr);
\draw(x3)[fill=white] circle(\vr);
\draw(x4)[fill=white] circle(\vr);
\draw(x5)[fill=white] circle(\vr);
\draw(x6)[fill=white] circle(\vr);
\end{scope}

\begin{scope}[xshift=6cm, yshift=0cm] % C4 4
\coordinate(x1) at (0.5,0.75);
\coordinate(x2) at (1,0);
\coordinate(x3) at (2,0);
\coordinate(x4) at (2.5,0.75);
\coordinate(x5) at (2,1.5);
\coordinate(x6) at (1,1.5);
% \edges		
\draw (x1) -- (x2) -- (x3) -- (x4) --(x5)-- (x6) --(x1);
%  vertices
\draw(x1)[fill=white] circle(\vr);
\draw(x2)[fill=white] circle(\vr);
\draw(x3)[fill=white] circle(\vr);
\draw(x4)[fill=white] circle(\vr);
\draw(x5)[fill=white] circle(\vr);
\draw(x6)[fill=white] circle(\vr);
\end{scope}

\begin{scope}[xshift=8cm, yshift=0cm] % C4 5
\coordinate(x1) at (0.5,0.75);
\coordinate(x2) at (1,0);
\coordinate(x3) at (2,0);
\coordinate(x4) at (2.5,0.75);
\coordinate(x5) at (2,1.5);
\coordinate(x6) at (1,1.5);
% \edges		
\draw (x1) -- (x2) -- (x3) -- (x4) --(x5)-- (x6) --(x1);
%  vertices
\draw(x1)[fill=white] circle(\vr);
\draw(x2)[fill=white] circle(\vr);
\draw(x3)[fill=white] circle(\vr);
\draw(x4)[fill=white] circle(\vr);
\draw(x5)[fill=white] circle(\vr);
\draw(x6)[fill=white] circle(\vr);
\end{scope}

\end{tikzpicture}
\caption{The graph $C_6^5$.}
	\label{fig:C_{6}^5}
\end{center}
\end{figure}

We conclude this section by demonstrating how colinear sets can be used while determining the general position number of graphs containing bridges.

\begin{proposition}
\label{prop:bridge}
Let $e=u_1u_2$ be a bridge of a connected graph $G$, where $\deg_G(u_1)\ge 2$ and $\deg_G(u_2)\ge 2$, and let $G_1$, $G_2$ the two components of $G-e$, where $u_i\in V(G_i)$, $i\in [2]$. Then
$\gp(G) \ge \xi_{G_1}(u_1) + \xi_{G_2}(u_2)$. Moreover, the equality holds in block graphs containing at least one bridge.
\end{proposition}

\proof
Let $X_i$ be a $u_i$-colinear set of $G_i$ with cardinality $\xi_{G_i}(u_i)$ for $i\in [2]$. Then it is straightforward to check that $X_1\cup X_2$ is a general position set of $G$, hence $\gp(G)\geq \xi_{G_1}(u_1) + \xi_{G_2}(u_2)$.

Consider the block graph $G$ with bridges, and let $e=u_1u_2$ be an arbitrary bridge of $G$. Let $G_1$ and $G_2$ be the two components of $G-e$, where $u_i\in V(G_i)$ and $i\in [2]$. From the above, we have $\gp(G)\geq \xi_{G_1}(u_1) + \xi_{G_2}(u_2)$. Let $S$ be the set of simplicial vertices of $G$ and set $S_i = S\cap V(G_i)$ for $i\in [2]$. Since the vertices from $S_i$, $i\in [2]$, are as described, we infer that $S_i$ is a $u_i$-colinear set of $G_i$. Clearly, $|S_1|+ |S_2| = |S|$ and by Theorem~\ref{thm:block-graphs} we know that $\gp(G) = |S|$. Putting all these fact together we have
$$\gp(G)\geq \xi_{G_1}(u_1) + \xi_{G_2}(u_2)\geq |S_1|+ |S_2| = |S| = \gp(G)\,,$$
which proves that $\gp(G) = \xi_{G_1}(u_1) + \xi_{G_2}(u_2)$.
\qed

Next, we give two examples of non-block graphs showing that the inequality of  Proposition~\ref{prop:bridge} can be sharp or not. Consider the graph $H$ as shown in Fig.~\ref{fig:H and H'}. It is straightforward to see that $\gp(H)=6$ and $\xi_{H_1}(u_1) + \xi_{H_2}(u_2) = 3+2 = 5$, hence $\gp(G)$ can be larger than $\xi_{G_1}(u_1) + \xi_{G_2}(u_2)$. On the other hand, to see that $\gp(G)$ can be equal to $\xi_{G_1}(u_1) + \xi_{G_2}(u_2)$ on non-block graphs with bridges, consider the graph $H'$ as shown in Fig.~\ref{fig:H and H'}. We have $\gp(H')=6$ and $\xi_{H'_1}(u_1) + \xi_{H'_2}(u_2) = 3+3 = 6$.

\begin{figure}[ht!]
\begin{center}
\begin{tikzpicture}[scale=1.0,style=thick,x=1cm,y=1cm]
\def\vr{2pt}

\begin{scope}[xshift=0cm, yshift=0cm] % C4 1
\coordinate(x1) at (0,0);
\coordinate(x2) at (1,0);
\coordinate(x3) at (2,0);
\coordinate(x4) at (3,0);
\coordinate(x5) at (4,0);
\coordinate(x6) at (5,0);
\coordinate(x7) at (7,0);
\coordinate(x8) at (6,0.75);
\coordinate(x9) at (3,0.75);
\coordinate(x10) at (1,0.75);
\coordinate(x11) at (1,-0.75);
\coordinate(x12) at (3,-0.75);
\coordinate(x13) at (6,-0.75);

% \edges		
\draw (x1) -- (x2) -- (x3) -- (x4) --(x5)-- (x6) --(x8) --(x7)--(x13)--(x6);
\draw (x1) -- (x10) -- (x3) -- (x11)-- (x1);
\draw (x3) -- (x9) -- (x5) -- (x12)-- (x3);

%  vertices
\draw(x1)[fill=white] circle(\vr);
\draw(x2)[fill=white] circle(\vr);
\draw(x3)[fill=white] circle(\vr);
\draw(x4)[fill=white] circle(\vr);
\draw(x5)[fill=white] circle(\vr);
\draw(x6)[fill=white] circle(\vr);
\draw(x7)[fill=white] circle(\vr);
\draw(x8)[fill=white] circle(\vr);
\draw(x9)[fill=white] circle(\vr);
\draw(x10)[fill=white] circle(\vr);
\draw(x11)[fill=white] circle(\vr);
\draw(x12)[fill=white] circle(\vr);
\draw(x13)[fill=white] circle(\vr);
% text
\node at (4.1,0.3) {$u_1$};
\node at (4.9,0.3) {$u_2$};
\node at (3,-1.3) {$H$};
\end{scope}

\begin{scope}[xshift=7.8cm, yshift=0cm] % C4 1
\coordinate(x1) at (0,0);
\coordinate(x2) at (0.75,0.75);
\coordinate(x3) at (1.75,0.75);
\coordinate(x4) at (2.5,0);
\coordinate(x5) at (1.75,-0.75);
\coordinate(x6) at (0.75,-0.75);
\coordinate(x7) at (3.5,0);
\coordinate(x8) at (4.5,0);
\coordinate(x9) at (5.25,-0.75);
\coordinate(x10) at (6,0);
\coordinate(x11) at (5.25,0.75);

% \edges		
\draw (x1) -- (x2) -- (x3) -- (x4) --(x5)-- (x6) --(x1) --(x3)--(x5)--(x1)--(x4)--(x6)--(x3);
\draw (x2)--(x4) -- (x7) -- (x8)--(x9)--(x10)--(x11)--(x8)--(x10);
\draw (x5) -- (x2)--(x6);
\draw (x9) -- (x11);
\draw (x3) -- (x7)--(x5);

%  vertices
\draw(x1)[fill=white] circle(\vr);
\draw(x2)[fill=white] circle(\vr);
\draw(x3)[fill=white] circle(\vr);
\draw(x4)[fill=white] circle(\vr);
\draw(x5)[fill=white] circle(\vr);
\draw(x6)[fill=white] circle(\vr);
\draw(x7)[fill=white] circle(\vr);
\draw(x8)[fill=white] circle(\vr);
\draw(x9)[fill=white] circle(\vr);
\draw(x10)[fill=white] circle(\vr);
\draw(x11)[fill=white] circle(\vr);
% text
\node at (3.5,0.3) {$u_1$};
\node at (4.4,0.3) {$u_2$};
\node at (4,-1.3) {$H'$};
\end{scope}

\end{tikzpicture}
\caption{The graphs $H$ and $H'$.}
	\label{fig:H and H'}
\end{center}
\end{figure}

%%%%%%%%%%%%%%%%%%%%%%
\section{Sierpi\'nski products of arbitrary graphs}
\label{sec:arbitrary}
%%%%%%%%%%%%%%%%%%%%%%

In this section we bound the (lower) Sierpi\'nski general position number for general graphs. Then we give formulas for the (lower) Sierpi\'nski general position number for the \SP\ of graphs with $K_2$ as the first factor. We begin with the following simple, but useful lemma.

\begin{lemma}
\label{lem:useful}
Let $X$ be a general position set of $K_2\otimes_f H$, where $V(K_2) = [2]$ and $f(2) = u\in V(H)$. If $|X\cap V(1H)| \ge 2$ and $(1,u)\in X$, then $X\cap V(2H) = \emptyset$.
\end{lemma}

\proof
Let $(1,v)$ be a vertex of $X$ different from $(1,u)$. Then every vertex from $2H$, say $(2,x)$, is an endvertex of a shortest $(1,v),(2,x)$-path passing through $(1,u)$, hence $(2,x)\notin X$.
\qed

\begin{theorem}
\label{thm:S-bounds}
If $G$ and $H$ are two connected graphs of order at least $2$, then
\[
\gp(H) \leq \gpsl(G,H)\leq \gps(G,H)\leq n(G)\gp(H).
\]
Moreover, $\gps(G,H) = n(G)\gp(H)$ if and only if $\gp(H) = \xi(H)$.
\end{theorem}

\proof
Let $G \otimes _f H$ be an arbitrary \SP\ of graphs $G$ and $H$. By Theorem~\ref{thm:convex}, the graph $G \otimes _f H$ contains $n(G)$ convex subgraphs $gH$, $g\in V(G)$. Since a convex subgraph is an  isometric subgraph, these subgraphs form an isometric cover. Therefore, by Theorem~\ref{thm:iso-cover}, we have $\gps(G,H)\leq n(G)\gp(H)$, hence the rightmost inequality holds. The convexity of the subgraphs $gH$ also implies the left inequality, while the middle inequality is obvious.

Assume now that $\gps(G,H) = n(G)\gp(H)$. Consider an arbitrary function $f\in H^G$. Using once more the convexity of the subgraph $gH$, there exists a gp-set $X$ of $G\otimes_f H$, such that $|V(gH) \cap X| = \gp(H)$ for each $g\in V(G)$. Consider a fixed subgraph $gH$ and let $gg'\in E(G)$. Since
$|V(gH) \cap X| = \gp(H)$ and $\gps(G,H) = n(G)\gp(H)$, Lemma~\ref{lem:useful} implies that $(g,f(g'))\notin X$. Setting $u = (g,f(g'))$ we obtain, having in mind Lemma~\ref{lem:useful} again, that $\xi_{gH}(u) = \gp(H)$ which in turn implies that $\xi(H) = \gp(H)$.

Assume next that $\xi(H) = \gp(H)$. Let $u$ be a vertex of $H$ such that $\xi_H(u) = \xi(H)$ and let $f\in H^G$ be the constant function $f(g)=u$ for any $g\in V(G)$. Each $gH$ has at least $\xi(H)$ vertices from some gp-set of the graph $G \otimes _f H$, hence $\gp(G \otimes _f H)\geq n(G)\xi(H) = n(G)\gp(H)$. Thus we have $\gps(G,H)\geq n(G)\gp(H)$.
Since each $gH$ is convex in $G \otimes _f H$, Theorem~\ref{thm:iso-cover} gives $\gp(G \otimes _f H)\leq n(G)\gp(H)$ and then $\gps(G,H)\leq n(G)\gp(H)$.
From the above, we conclude that $\gps(G,H) = n(G)\gp(H)$.
\qed

\begin{corollary}
If $G$ is a connected graph with $n(G)\ge 2$ and $T$ is a tree with $n(T)\ge 3$, then $\gps(G,T) = n(G)\ell(T)$.
\end{corollary}

\proof
Since $n(T)\ge 3$, Proposition~\ref{prop:xi-in-tree} gives $\gp(T) = \ell(T) = \xi(T)$. Hence Theorem~\ref{thm:S-bounds} yields the conclusion.
\qed

\begin{theorem}
\label{thm:K2-H}
If $H$ is a connected graph with $n(H)\ge 2$, then the following assertions hold.

(i) $\gpsl(K_2,H) = 2\xi^{-}(H)$.

(ii) $\gps(K_2,H) = 2\xi(H)$.
\end{theorem}

\proof
Set $V(K_2) = [2]$ for the rest of the proof.

(i) Let $f\in H^{K_2}$ be an arbitrary function, say $f(1) = u$ and $f(2) = v$. By Proposition~\ref{prop:bridge} we have $\gp(K_2 \otimes _f H) \ge \xi_{H}(v) + \xi_{H}(u) \ge \xi^{-}(H) + \xi^{-}(H)$. It follows that $\gpsl(K_2,H) \ge 2\xi^{-}(H)$.

To prove the reverse inequality, let $u$ be a vertex of $H$ such that $\xi_H(u) = \xi^{-}(H)$. Define $f\in H^{K_2}$ by $f(1) = f(2) = u$ and let $X$ be a gp-set of $K_2 \otimes _f H$. If $V(1H)\cap X = \emptyset$ or $V(2H)\cap X = \emptyset$, then $\gp(K_2 \otimes _f H)\leq \gp(H)$ holds and then Theorem~\ref{thm:2xi>gp} gives $\gp(K_2 \otimes _f H)\leq 2\xi^{-}(H)$. Assume second that $V(1H)\cap X \ne \emptyset$ and $V(2H)\cap X \ne \emptyset$. Then we claim that $|V(1H)\cap X|\le \xi_H(u)$ and $|V(2H)\cap X|\le \xi_H(u)$. By establishing this claim, $\gp(K_2 \otimes _f H)\leq 2\xi^{-}(H)$ follows. Suppose on the contrary that this is not the case and assume without loss of generality that $|V(1H)\cap X|\ge \xi_H(u)+1$. Then there exist two vertices $(1,h), (1,h')\in V(1H)\cap X$, where $h\ne h'$, such that the vertices $(1,h)$, $(1,h')$, and $(2,u)$ lie on a common shortest path.  If $(2, h'')$ is an arbitrary vertex from $V(2H)\cap X$, then the vertices $(1,h)$, $(1,h')$, and $(2, h'')$ lie on a common shortest path. It is impossible because we have assumed that $V(2H)\cap X \ne \emptyset$. This contradiction proves the claim.

(ii) Let $f\in H^{K_2}$ be an arbitrary function, say $f(1) = u$ and $f(2) = v$ and let $X$ be an arbitrary gp-set of $K_2 \otimes _f H$.
We claim that $\gp(K_2 \otimes _f H) \leq 2\xi(H)$. Suppose on the contrary that $\gp(K_2 \otimes _f H) \geq 2\xi(H) + 1$. Since the order of $H$ is at least $2$, Theorem~\ref{thm:2xi>gp} implies that $2\xi(H)\geq \gp(H)\geq 2$. Then $\gp(K_2 \otimes _f H) \geq \gp(H) + 1 \geq 3$ and thus
$|X\cap V(1H)|\geq 2$ or $|X\cap V(2H)|\geq 2$.

Assume, without loss of generality, that $|X\cap V(1H)|\geq 2$. If $(1,v)\in X$, then Lemma~\ref{lem:useful} implies that $X\cap V(2H) = \emptyset$ and then, using Theorem~\ref{thm:2xi>gp} once more, $\gp(K_2 \otimes _f H) \leq \gp(H) \le 2\xi^{-}(H) \le 2\xi(H)$, a contradiction. Hence assume in the rest that $(1,v)\notin X$. Since $\gp(K_2 \otimes _f H) \geq \gp(H) + 1$, we get $|X\cap V(2H)|\geq 2$. But as argued in (i), this is not possible for $(2,u)\in X$. Hence $|X\cap V(2H)|\leq \xi_H(u)\leq \xi(H)$. By symmetry we also have that if $|X\cap V(2H)|\geq 2$, then $|X\cap V(1H)|\leq \xi_H(v)\leq \xi(H)$ holds. We can thus conclude that $\gps(K_2,H)\leq 2\xi(H)$.

To complete the argument we need to demonstrate that $\gps(K_2,H)\geq 2\xi(H)$. Let $x$ be a vertex of $H$ such that $\xi_H(x) = \xi(H)$ and let the function $f\in H^{K_2}$ be defined by $f(1) = f(2) = x$. By Proposition~\ref{prop:bridge}, it follows that $\gp(K_2 \otimes _f H)\geq 2\xi_H(x) = 2\xi(H)$.
Then we have $\gps(K_2,H) \geq 2\xi(H)$, and we are done.
\qed

%%%%%%%%%%%%%%%%%%%%%%
\section{Sierpi\'nski products of complete graphs}
\label{sec:complete graphs}
%%%%%%%%%%%%%%%%%%%%%%

In this section we consider the \SP\ of complete graphs $K_m \otimes _f K_n$. In the first main result we determine their Sierpi\'nski general position number. For the lower Sierpi\'nski general position number we prove that $\gpsl(K_m,K_n) = m(n-m+1)$ provided that $n\ge 2m-2$. The latter assumption is required as demonstrated by our last result which asserts that $\gpsl(K_6,K_9)=25$.

Throughout this section we will assume that $V(K_m) = [m]$ and $V(K_n) = \{x_i:\ i\in[n]\}$.

\begin{theorem}
\label{thm:gps}
If $m,n\geq 2$, then
\[
\gps(K_m,K_n) = m(n-1)\,.
\]
\end{theorem}

\proof
Let $f\in K_n^{K_m}$ be defined by $f(i)=x_1$ for $i\in [m]$. Let
$X = \{(i,x_j):\ i\in[m],\ 2\leq j\leq n\}$. Since $f$ is a constant function, for any two vertices $(i,x_j), (i',x_j')\in X$, where $i\neq i'$, we have $d_{K_m \otimes _f K_n}((i,x_j), (i',x_j')) = 3$. By Theorem~\ref{thm:gpsets}, $X$ is a general position set, hence $\gp(K_m \otimes _f K_n)\geq  m(n-1)$. It follows that $\gps(K_m,K_n)\geq \gp(K_m \otimes _f K_n) \geq m(n-1)$.

To prove equality, suppose on the contrary that $\gps(K_m,K_n) \geq m(n-1) + 1$. It means that there exists a function $f'\in K_n^{K_m}$ and an index $i\in [m]$,  such that $V(iK_n)$ completely lies in some general position set $S$ of $K_m \otimes _{f'} K_n$. Since $m, n\geq 2$, there exists a vertex $(i',x_j)\in S$, where $i\neq i'$. Assume that $f'(i') = x_p$ and let $q\in [n]$ be such that $q\ne p$. Then the vertices $(i,x_q)\in S$, $(i,x_p)\in S$, and $(i',x_j)\in S$ lie on a common shortest path (of length $2$) of $K_m \otimes _{f'} K_n$, a contradiction. This contradiction implies that $\gp(K_m \otimes _{f'} K_n) \leq m(n-1)$, hence we have $\gps(K_m,K_n) = m(n-1)$.
\qed

In the remainder we focus on the lower Sierpi\'nski general position number.

\begin{lemma}
\label{lem:n-1-vertex}
Let $m,n\geq 2$ and let $f\in K_n^{K_m}$. Then the following hold.

(i) If $m \le n$, then $\gp(K_m \otimes _f K_n) \ge v_{n-1}(K_m \otimes _f K_n) \ge (n-m+1)m$. In particular,
$$\gpsl(K_m,K_n)\ge (n-m+1)m\,.$$

(ii) If $m\ge n$, then $\gp(K_m \otimes _f K_n) \ge \max \{2(n-1), m\}$. In particular,
$$\gpsl(K_m,K_n)\ge \max \{2(n-1), m\}\,.$$

\end{lemma}

\proof
Let $f\in K_n^{K_m}$, and set $Z = K_m \otimes _f K_n$ for the rest of the proof.

(i) Note that a vertex of $Z$ is of degree $n-1$ if and only if it is not incident to a connecting edge. For each $i\in [m]$, let $X_i$ be the set of vertices from $V(Z) \cap V(iK_n)$ of degree $n-1$. (Note that it is possible that $X_i = \emptyset$.) Then $X_i$ clearly induces a complete subgraph. Moreover, if $i\ne j$, and if $(i,x_k)\in X_i$ and $(j,x_\ell)\in X_j$, then $d_{Z}((i,x_k), (j,x_\ell)) = 3$. Theorem~\ref{thm:gpsets} implies that $\bigcup_{i=1}^m X_i$ is a general position set. This in turn yields $\gp(Z) \ge v_{n-1}(Z)$ because each vertex not in  $\bigcup_{i=1}^m X_i$ is of degree at least $n$. Moreover, since be definition, each $iK_n$ contains at most $m-1$ vertices incident to connecting edges, we have $|X_i|\ge n-m+1$, so we also have $v_{n-1}(Z) \ge (n-m+1)m$.

The above argument holds true for any function $f\in K_n^{K_m}$, hence we have $\gpsl(K_m,K_n)\ge (n-m+1)m$.

(ii) Assume first that $f$ is a constant function, without loss of generality let $f(i) = x_1$ for $i\in [m]$. By Theorem~\ref{thm:gpsets} we then get that $V(Z)\setminus \{(1,x_1), \ldots, (m,x_1)\}$ is a general position set, hence in this case $\gp(Z) \ge 2(n-1)$. Assume second that $f$ is not an identity function and assume without loss of generality that $f(1) = x_1$ and $f(2) = x_2$. Now we claim that $S = V(1K_n) \cup V(2K_n) \setminus \{(1,x_2), (2,x_1)\}$ is a general position set. To prove this claim it suffices to show (again in view of Theorem~\ref{thm:gpsets}) that if $(1,x_i), (2,x_j)\in S$, then $d_{Z}((1,x_i), (2,x_j)) = 3$. Suppose on the contrary that $d_{Z}((1,x_i), (2,x_j)) = 2$. Then there exists a vertex $(p,x_k)$, where $k\ne 1,2$, such that
$(1,x_i)(p,x_k)\in E(Z)$ and $(2,x_j)(p,x_k)\in E(Z)$. But this means that $f(1) = x_k$ and $f(2) = x_k$, a contradiction since we have assumed that $f(1)\ne f(2)$. This contradiction proves the claim. We conclude that $S$ is a general position set of $Z$ and as $|S| = 2(n-1)$ we have proved that $\gp(Z) \ge 2(n-1)$.

To prove that $\gp(Z) \ge m$, we may assume without loss of generality that the function $f$ is non-decreasing with respect to the indices of $x_i$s. That is, let $k_1\ge \cdots \ge k_n \ge 0$, where $k_1+\cdots + k_n = m$, and set $k_0 = 0$. Then for $i\in [n]$,
$$f(j) = x_i\quad {\rm for}\quad j\in \{k_1+\cdots + k_{i-1}+1, \ldots, k_1+\cdots + k_{i}\}\,.$$
For $i\in [n]$, set
$$X_i = \big\{(j,x_i):\ j\in \{k_1+\cdots + k_{i-1}+1, \ldots, k_1+\cdots + k_{i}\}\big\}\,.$$
Note that for some $i$ we can have $X_i = \emptyset$, but in any case $\sum_{i=1}^n |X_i| = m$.

We claim that $X = \bigcup_{i=1}^nX_i$ is a general position set of $Z$. Note first that $X_i$, $i\in [n]$, induces a complete subgraph $K_{k_i}$ of $Z$. Therefore, in view of Theorem~\ref{thm:gpsets}, to prove that $X$ is a general position set, it suffices to demonstrate that if $(j,x_i)\in X_i$ and $(j',x_{i'})\in X_{i'}$, where $i\ne i'$, then $d_Z((j,x_i), (j',x_{i'})) = 3$. Since $(j,x_i)$ and $(j',x_{i'})$ are not adjacent, suppose on the contrary that $d_Z((j,x_i), (j',x_{i'})) = 2$. Then there exists a vertex $(p,x_q)$, where $p\ne j, j'$, such that $(j,x_i)(p,x_q)\in E(Z)$ and $(j',x_{i'})(p,x_q)\in E(Z)$. Then it follows that $f(j)=x_q$ and $f(j')=x_q$. But this is not possible, since $(j,x_i)\in X_i$ and $(j',x_{i'})\in X_{i'}$, where $i\ne i'$, which in particular means that $f(j)\ne f(j')$. This contradiction proves the claim, hence $\gp(Z) \ge m$.

Again, the above arguments hold true for any function $f\in K_n^{K_m}$, hence we have $\gpsl(K_m,K_n)\ge \max \{2(n-1), m\}$.
\qed

For the lower Sierpi\'nski general position number of two complete graphs we have the following result, where we need to assume that the second factor is relatively large with respect to the first factor.

\begin{theorem}
\label{thm:gpsl}
If $m\geq 2$ and $n\ge 2m-2$, then
\[
\gpsl(K_m,K_n) = m(n-m+1)\,.
\]
\end{theorem}

\proof
Using Lemma~\ref{lem:n-1-vertex}(i) once more we have $\gpsl(K_m,K_n) \geq m(n-m+1)$.

To prove that $\gpsl(K_m,K_n) = m(n-m+1)$, consider the function $f\in K_n^{K_m}$ defined by $f(i)=x_i$ for $i\in [m]$. Setting $G = K_m\otimes_{f} K_n$ we are going to show that $\gp(G)\leq m(n-m+1)$.

Suppose on the contrary that $\gp(G)\geq m(n-m+1) + 1$ and let $R$ be a gp-set of $G$. Then it follows that there exists an index $i\in[m]$ such that $|V(iK_n)\cap R|\geq n-m+2\ge 2$. Assume, without loss of generality, that $|V(1K_n)\cap R|$ is as large as possible.

We claim that if $(1,x_j)\in V(1K_n)\cap R$, then $V(jK_n)\cap R=\emptyset$, $j\in [m]$. Indeed, since $|V(1K_n)\cap R|\geq 2$, there exists a vertex $(1, x_{j'})\in R$, where $j'\ne j$ and $j'\in[m]$. Then the vertices $(1, x_{j'})$, $(1, x_{j})$, $(j, x_1)$, and $(j,x_p)\in V(jK_n)$ form a shortest path, hence we conclude that $(j,x_p)\notin R$ for $p\in [n]$.

As we have assumed that $\gp(G)\geq m(n-m+1) + 1$, we can write $|V(1K_n)\cap R| = n-m+1+k$, where $k\ge 1$ (and $k\le m-1$). Then $|V(sK_n)\cap R|\leq n-m+1+k$, where $2\leq s\leq m$. By the above, $G$ contains at least $k$ copies $iK_n$ such that $V(iK_n)\cap R = \emptyset$. Then we have
\begin{align*}
|R|&\leq (n-m+1+k)+(m-k-1)(n-m+1+k)\\
& = m(n-m+1)+k(2m-n-1-k)\\
&< m(n-m+1) + 1.
\end{align*}
Here the last inequality holds because $n\ge 2m - 2$ and $k\ge 1$ which implies that $2m-n-1-k \le 2 - 1 - k = 1 - k \le 0$. This contradiction implies that $\gp(G)\leq m(n-m+1)$. We conclude that $\gpsl(K_m,K_n) = m(n-m+1)$.
\qed

In Theorem~\ref{thm:gpsl} we have assumed that $n\ge 2m-2$. The following result explains why this assumption cannot be avoided in general.

\begin{proposition}
\label{prop:K6K9}
$\gpsl(K_6,K_9)=25$.
\end{proposition}

\proof
By Lemma~\ref{lem:n-1-vertex}(i), we have $\gpsl(K_6,K_9)\geq 24$. To prove that $\gpsl(K_6,K_9) = 25$, we consider two cases depending on the function $f\in K_9^{K_6}$.

\medskip\noindent
{\bf Case 1}. $f$ is injective. \\
In this case we may assume, without loss of generality, that $f(i)=x_i$ for $i\in[6]$. The obtained graph $K_6\otimes_f K_9$ is shown in Fig.~\ref{fig:K_6-K_9}, where the edges of complete subgraphs $iK_9$, $i\in [6]$, and the edges of the factors $K_6$ and $K_9$ are not drawn for a clearer picture.

\begin{figure}[ht!]
\begin{center}
\begin{tikzpicture}[scale=0.8,style=thick,x=1.8cm,y=1cm]
\def\vr{3pt}

\begin{scope}[xshift=0cm, yshift=0cm] % C4 1
\coordinate(x1) at (1,0);
\coordinate(x2) at (2,0);
\coordinate(x3) at (3,0);
\coordinate(x4) at (4,0);
\coordinate(x5) at (5,0);
\coordinate(x6) at (6,0);
\coordinate(u1) at (1,-2);
\coordinate(u2) at (2,-2);
\coordinate(u3) at (3,-2);
\coordinate(u4) at (4,-2);
\coordinate(u5) at (5,-2);
\coordinate(u6) at (6,-2);
\coordinate(y1) at (-0.5,0);
\coordinate(y2) at (-0.5,0.8);
\coordinate(y3) at (-0.5,1.6);
\coordinate(y4) at (-0.5,2.4);
\coordinate(y5) at (-0.5,3.2);
\coordinate(y6) at (-0.5,4.0);
\coordinate(y7) at (-0.5,4.8);
\coordinate(y8) at (-0.5,5.6);
\coordinate(y9) at (-0.5,6.4);

% \edges		
%\draw (x1) -- (x2) -- (x3) -- (x4) --(x1);
\draw (-0.5,-1.5) -- (7,-1.5);
\draw (0,-2.5) -- (0,7.5);

%  vertices
\draw(x1)[fill=black] circle(\vr);
\draw(x2)[fill=white] circle(\vr);
\draw(x3)[fill=white] circle(\vr);
\draw(x4)[fill=white] circle(\vr);
\draw(x5)[fill=white] circle(\vr);
\draw(x6)[fill=white] circle(\vr);
\draw(u1)[fill=white] circle(\vr);
\draw(u2)[fill=white] circle(\vr);
\draw(u3)[fill=white] circle(\vr);
\draw(u4)[fill=white] circle(\vr);
\draw(u5)[fill=white] circle(\vr);
\draw(u6)[fill=white] circle(\vr);
\draw(y1)[fill=white] circle(\vr);
\draw(y2)[fill=white] circle(\vr);
\draw(y3)[fill=white] circle(\vr);
\draw(y4)[fill=white] circle(\vr);
\draw(y5)[fill=white] circle(\vr);
\draw(y6)[fill=white] circle(\vr);
\draw(y7)[fill=white] circle(\vr);
\draw(y8)[fill=white] circle(\vr);
\draw(y9)[fill=white] circle(\vr);

%\draw (x1) .. controls (-1.8,6.35).. (x6);
\draw (1,3.2) ellipse (0.3 and 3.6);
\node at (1,-0.8) {$1K_9$};
\end{scope}

\begin{scope}[xshift=0cm, yshift=0.8cm] % C4 2
\coordinate(x1) at (1,0);
\coordinate(x2) at (2,0);
\coordinate(x3) at (3,0);
\coordinate(x4) at (4,0);
\coordinate(x5) at (5,0);
\coordinate(x6) at (6,0);
% \edges		
%\draw (x1) -- (x2) -- (x3) -- (x4) --(x1);
%  vertices
\draw(x1)[fill=white] circle(\vr);
\draw(x2)[fill=black] circle(\vr);
\draw(x3)[fill=white] circle(\vr);
\draw(x4)[fill=white] circle(\vr);
\draw(x5)[fill=white] circle(\vr);
\draw(x6)[fill=white] circle(\vr);

\draw (2,2.4) ellipse (0.3 and 3.6);
\node at (1, -3.5) {$1$};
\node at (2, -3.5) {$2$};
\node at (3, -3.5) {$3$};
\node at (4, -3.5) {$4$};
\node at (5, -3.5) {$5$};
\node at (6, -3.5) {$6$};
\node at (-0.8, -0.8) {$x_1$};
\node at (-0.8, 0) {$x_2$};
\node at (-0.8, 0.8) {$x_3$};
\node at (-0.8, 1.6) {$x_4$};
\node at (-0.8, 2.4) {$x_5$};
\node at (-0.8, 3.2) {$x_6$};
\node at (-0.8, 4.0) {$x_7$};
\node at (-0.8, 4.8) {$x_8$};
\node at (-0.8, 5.6) {$x_9$};
\node at (2,-1.6) {$2K_9$};
\end{scope}

\begin{scope}[xshift=0cm, yshift=1.6cm] % C4 3
\coordinate(x1) at (1,0);
\coordinate(x2) at (2,0);
\coordinate(x3) at (3,0);
\coordinate(x4) at (4,0);
\coordinate(x5) at (5,0);
\coordinate(x6) at (6,0);

% \edges		
%\draw (x1) -- (x2) -- (x3) -- (x4) --(x1);
%  vertices
\draw(x1)[fill=white] circle(\vr);
\draw(x2)[fill=white] circle(\vr);
\draw(x3)[fill=black] circle(\vr);
\draw(x4)[fill=white] circle(\vr);
\draw(x5)[fill=white] circle(\vr);
\draw(x6)[fill=white] circle(\vr);

\draw (3,1.6) ellipse (0.3 and 3.6);
\node at (3,-2.4) {$3K_9$};
\end{scope}

\begin{scope}[xshift=0cm, yshift=2.4cm] % C4 4
\coordinate(x1) at (1,0);
\coordinate(x2) at (2,0);
\coordinate(x3) at (3,0);
\coordinate(x4) at (4,0);
\coordinate(x5) at (5,0);
\coordinate(x6) at (6,0);
% \edges		
%\draw (x1) -- (x2) -- (x3) -- (x4) --(x1);
%  vertices
\draw(x1)[fill=white] circle(\vr);
\draw(x2)[fill=white] circle(\vr);
\draw(x3)[fill=white] circle(\vr);
\draw(x4)[fill=black] circle(\vr);
\draw(x5)[fill=white] circle(\vr);
\draw(x6)[fill=white] circle(\vr);

\draw (4,0.8) ellipse (0.3 and 3.6);
\node at (4,-3.2) {$4K_9$};
\end{scope}

\begin{scope}[xshift=0cm, yshift=3.2cm] % C4 5
\coordinate(x1) at (1,0);
\coordinate(x2) at (2,0);
\coordinate(x3) at (3,0);
\coordinate(x4) at (4,0);
\coordinate(x5) at (5,0);
\coordinate(x6) at (6,0);
% \edges		
%\draw (x1) -- (x2) -- (x3) -- (x4) --(x1);
%  vertices
\draw(x1)[fill=white] circle(\vr);
\draw(x2)[fill=white] circle(\vr);
\draw(x3)[fill=white] circle(\vr);
\draw(x4)[fill=white] circle(\vr);
\draw(x5)[fill=black] circle(\vr);
\draw(x6)[fill=white] circle(\vr);

\draw (5,0) ellipse (0.3 and 3.6);
\node at (5,-4) {$5K_9$};
\end{scope}

\begin{scope}[xshift=0cm, yshift=4cm] % C4 6
\coordinate(x1) at (1,0);
\coordinate(x2) at (2,0);
\coordinate(x3) at (3,0);
\coordinate(x4) at (4,0);
\coordinate(x5) at (5,0);
\coordinate(x6) at (6,0);
% \edges		
%\draw (x1) -- (x2) -- (x3) -- (x4) --(x1);
%  vertices
\draw(x1)[fill=black] circle(\vr);
\draw(x2)[fill=black] circle(\vr);
\draw(x3)[fill=black] circle(\vr);
\draw(x4)[fill=black] circle(\vr);
\draw(x5)[fill=black] circle(\vr);
\draw(x6)[fill=white] circle(\vr);

\draw (6,-0.8) ellipse (0.3 and 3.6);
\node at (6,-4.8) {$6K_9$};
\end{scope}

\begin{scope}[xshift=0cm, yshift=4.8cm] % C4 7
\coordinate(x1) at (1,0);
\coordinate(x2) at (2,0);
\coordinate(x3) at (3,0);
\coordinate(x4) at (4,0);
\coordinate(x5) at (5,0);
\coordinate(x6) at (6,0);
% \edges		
%\draw (x1) -- (x2) -- (x3) -- (x4) --(x1);
%  vertices
\draw(x1)[fill=black] circle(\vr);
\draw(x2)[fill=black] circle(\vr);
\draw(x3)[fill=black] circle(\vr);
\draw(x4)[fill=black] circle(\vr);
\draw(x5)[fill=black] circle(\vr);
\draw(x6)[fill=white] circle(\vr);
\end{scope}

\begin{scope}[xshift=0cm, yshift=5.6cm] % C4 8
\coordinate(x1) at (1,0);
\coordinate(x2) at (2,0);
\coordinate(x3) at (3,0);
\coordinate(x4) at (4,0);
\coordinate(x5) at (5,0);
\coordinate(x6) at (6,0);
% \edges		
%\draw (x1) -- (x2) -- (x3) -- (x4) --(x1);
%  vertices
\draw(x1)[fill=black] circle(\vr);
\draw(x2)[fill=black] circle(\vr);
\draw(x3)[fill=black] circle(\vr);
\draw(x4)[fill=black] circle(\vr);
\draw(x5)[fill=black] circle(\vr);
\draw(x6)[fill=white] circle(\vr);
\end{scope}

\begin{scope}[xshift=0cm, yshift=6.4cm] % C4 9
\coordinate(x1) at (1,0);
\coordinate(x2) at (2,0);
\coordinate(x3) at (3,0);
\coordinate(x4) at (4,0);
\coordinate(x5) at (5,0);
\coordinate(x6) at (6,0);
% \edges		
%\draw (x1) -- (x2) -- (x3) -- (x4) --(x1);
%  vertices
\draw(x1)[fill=black] circle(\vr);
\draw(x2)[fill=black] circle(\vr);
\draw(x3)[fill=black] circle(\vr);
\draw(x4)[fill=black] circle(\vr);
\draw(x5)[fill=black] circle(\vr);
\draw(x6)[fill=white] circle(\vr);
\end{scope}

\draw (1.05,0.78) -- (1.96,0.1);
\draw (1.05,1.52) .. controls (2,0.5).. (2.95,0);
\draw (1.05,2.35) .. controls (2.45,0.75).. (3.95,0.05);
\draw (2.05,1.58) -- (2.96,0.85);
\draw (1.05,3.15) .. controls (2.6,1.4).. (4.95,0);
\draw (2.05,2.45) .. controls (3.1,1.85).. (3.95,0.88);
\draw (1.05,3.95) .. controls (3.3,1.65).. (5.95,0);
\draw (3.05,2.38) -- (3.95,1.65);
\draw (2.05,3.18) .. controls (3.5,2.4).. (4.95,0.88);
\draw (3.05,3.2) .. controls (4.1,2.65).. (4.95,1.62);
\draw (2.05,4) .. controls (4.1,3).. (5.95,0.88);
\draw (4.05,3.22) -- (4.95,2.4);
\draw (3.05,4) .. controls (4.35,3.4).. (5.95,1.62);
\draw (4.05,4.05) .. controls (5,3.5).. (5.95,2.4);
\draw (5.05,4.05) -- (5.95,3.2);

\end{tikzpicture}
\caption{The \SP\ graph $K_6\otimes_f K_9$, where  $f(i)=x_i$ for $i\in[6]$}
	\label{fig:K_6-K_9}
\end{center}
\end{figure}

We can check that the set consisting of black vertices from Fig.~\ref{fig:K_6-K_9} is a general position set of cardinality $25$ of $K_6\otimes_f K_9$. Hence $\gp(K_6\otimes_f K_9) \ge 25$. To prove the reverse inequality, suppose there exists a general position set $X$ of $K_6\otimes_f K_9$ with $|X|\ge 26$. Let $X_i = X\cap V(iK_9)$, $i\in [6]$. We may assume without loss of generality that $|X_1| \ge |X_i|$ for $i\in \{2,\ldots, 6\}$. Clearly, $|X_1| \ge 5$. If $|X_1| = 5$, then for some $i\in \{2,\ldots, 6\}$ we have $X_i = \emptyset$. This implies that  $|X|\le 25$. Hence we must have $|X_1| = 4 + k$, where $2\le k\le 5$. Then there are $k$ copies $iK_9$ with no vertex from $X$. It follows that $|X| \le (4 + k) + (6-1-k)(4+k) \le 25$. We can conclude that $\gp(K_6\otimes_f K_9) \le 25$ and thus $\gp(K_6\otimes_f K_9) = 25$.

\medskip\noindent
{\bf Case 2}. $f$ is not injective. \\
In this case we may assume, without loss of generality, that $f(1)=f(2)$. Then in each of the copies $iK_9$, $i\in \{3,4,5,6\}$, the graph $K_6\otimes_f K_9$ contains at least five vertices of degree $8$. Since in each of $1K_9$ and $2K_9$ there are at least four such vertices, we can see (using the argument from the proof of Lemma~\ref{lem:n-1-vertex}(i)) that $\gp(K_6\otimes_f K_9) \ge 8 + 20 = 28$.
\qed

%%%%%%%%%%%%%%%%%%%%%%%%%%%%
\section{Concluding remarks}
\label{S:remark}
%%%%%%%%%%%%%%%%%%%%%%%%%%%%

In Proposition~\ref{prop:bridge} we have bounded from below the general position number of graphs with bridges. It would be interesting to characterize the graphs that attain the equality in the proved bound.

In view of Theorem~\ref{thm:S-bounds} it would be interesting to characterize the graphs $G$ with $\gp(G) = \xi(G)$. Such graphs are, for instance, grid graphs $P_n \cp P_m$, $n,m \ge 3$, for which we have $\gp(P_n \cp P_m) = 4 = \xi(P_n \cp P_m)$. Similarly, in view of Theorem~\ref{thm:K2-H} it would be interesting to study the graphs $G$ with $\xi^{-}(G) = \xi(G)$.

In Theorem~\ref{thm:gpsl} we have proved that $\gpsl(K_m,K_n) = m(n-m+1)$ if $m\geq 2$ and $n\ge 2m-2$. It remains as an open problem to determine $\gpsl(K_m,K_n)$ for the cases when $n < 2m-2$. In Proposition~\ref{prop:K6K9} the particular case of $\gpsl(K_6,K_9)$ has been solved.

%%%%%%%%%%%%%%%%%%%%%%%%%%%%%%%%%%%%%%%%%%%%%%%%%%%%%%%%
\section*{Acknowledgements}
%%%%%%%%%%%%%%%%%%%%%%%%%%%%%%%%%%%%%%%%%%%%%%%%%%%%%%%%

We would like to express our sincere thanks to the three reviewers for their many helpful suggestions on how to improve the paper.

This work was supported by the Slovenian Research Agency ARIS (research core funding P1-0297 and projects N1-0285, N1-0355).

%%%%%%%%%%%%%%%%%%%%%%%%%%%%
\section*{Declaration of interests}
%%%%%%%%%%%%%%%%%%%%%%%%%%%%

The authors declare that they have no known competing financial interests or personal relationships that could have appeared to influence the work reported in this paper.

%%%%%%%%%%%%%%%%%%%%%%%%%%%%
\section*{Data availability}
%%%%%%%%%%%%%%%%%%%%%%%%%%%%

Our manuscript has no associated data.

\end{document}